\documentclass[11pt]{amsart}
\usepackage{amsfonts}
\usepackage{amsfonts,latexsym,rawfonts,amsmath,amssymb,amsthm,enumerate}
\usepackage[plainpages=false]{hyperref}

\usepackage{graphicx}

\RequirePackage{color}

\newcommand\ZZ{\mathbb{Z}}

\newcounter{marnote}

\def\Ric{{\rm Ric}}
\def\RR{{\mathbb{R}}}

\theoremstyle{plain}
  \newtheorem{theorem}{Theorem}
  \newtheorem{lemma}{Lemma}
  \newtheorem{proposition}{Proposition}[section]
  
  \newtheorem{corollary}{Corollary}[section]
  \newtheorem{definition}{Definition}

\theoremstyle{remark} 
  \newtheorem{remark}{Remark}

\renewcommand{\det}{\mbox{det}}
  
  \newcommand{\dist}{\mbox{dist}}

  \numberwithin{equation}{section}
  \numberwithin{figure}{section}

\begin{document}

\title[Nonlinear oblique operators]
{A degree theory for second order nonlinear elliptic operators with nonlinear oblique boundary conditions}

\date{}

\author[Y.Y. Li]{YanYan Li}

\address
	{School of Mathematical Sciences,
	Beijing Normal University,
	Beijing 100875, China}

\address{and}

\address{Department of Mathematics,
	Rutgers University,
	110 Frelinghuysen Road,
	Piscataway, NJ 08854-8019, USA.}

\email{yyli@math.rutgers.edu}

\author[J. Liu]{Jiakun Liu}

\address
	{Institute for Mathematics and its Applications,
	School of Mathematics and Applied Statistics,
	University of Wollongong,
	Wollongong, NSW 2522, Australia.}
	
\email{jiakunl@uow.edu.au}

\author[L. Nguyen]{Luc Nguyen}

\address
	{Mathematical Institute and St Edmund Hall,
	University of Oxford,
	Andrew Wiles Building,
	Radcliffe Observatory Quarter,
	Woodstock Road,
	Oxford OX2 6GG, UK.}
	
\email{luc.nguyen@maths.ox.ac.uk}

\thanks{ }
\thanks{\copyright 2012 by the authors. All rights reserved}

\subjclass[2000]{ }

\keywords{ }

\begin{abstract}
In this paper we introduce an integer-valued degree for second order fully nonlinear elliptic operators with nonlinear oblique boundary conditions. We also give some applications to the existence of solutions of certain nonlinear elliptic equations arising from a Yamabe problem with boundary and reflector problems. 
\end{abstract}

\dedicatory{Dedicated to Paul H. Rabinowitz on his 75th birthday with admiration}
\maketitle

\baselineskip=16.4pt
\parskip=3pt

%++++++++++++++++%

\section{Introduction}

Degree theories are very useful in the study of partial differential equations, for example, in the study of existence and multiplicities of solutions, eigenvalue and bifurcation problems. See for example \cite{CrandallRabinowitz, Krasno, NirBook, Rabinowitz1, Rabinowitz2}.

In \cite{Li}, the first named author introduced a degree theory for second order nonlinear elliptic operators with Dirichlet boundary conditions. It is natural to ask for a degree theory for other boundary operators. Problems with nonlinear oblique boundary conditions have been considered in the literature for some time, see e.g. \cite{Chen, JiangTruXiang, JinLiLi, LL, LLCounter, LieTru, LionsTruUrbas, LT, MaQiu, Urbas2, Urbas}.

For example, in the study of boundary Yamabe problems \cite{Chen, JinLiLi, LL, LLCounter}, one considers the boundary condition
  \begin{equation}\label{1001}
   h_{u^{\frac{4}{n-2}}g} := u^{-\frac{n}{n-2}}\left[\frac{\partial u}{\partial\nu}+\frac{n-2}{2}h_gu\right]=c \quad\mbox{on }\partial M,
  \end{equation}
where $\partial M$ is the boundary of a smooth Riemannian manifold $(M,g)$ of dimension $n\geq 3$, $\nu$ is the outer unit normal to $\partial M$ and $h_g$ is the mean curvature of $\partial M$. \eqref{1001} is a semi-linear Neumann boundary condition. 

More recently, in the study of a near field reflector problem \cite{LT} one has the boundary condition
  \begin{equation}\label{1002}
   T_u(\Omega) = \Omega^*, 
  \end{equation}
where $\Omega, \Omega^*$ are two bounded domains in $\mathbb{R}^n$, and $T_u$ is the reflection mapping given by
  \begin{equation*}
   T_u(x) = \frac{2Du}{|Du|^2-(u-Du\cdot x)^2},\quad x\in\Omega.
  \end{equation*}
The boundary condition \eqref{1002} is fully nonlinear, and in \cite{LT} it was shown that \eqref{1002} is oblique for any admissible solution $u$. 

The equations associated with \eqref{1001} and \eqref{1002} are Hessian and Monge-Amp\`ere types, respectively.

The main goal of the present paper is to define a degree theory, along the line of \cite{Li, LL} for fully nonlinear elliptic operators with fully nonlinear oblique boundary conditions. See Section \ref{Sec2} for the statement and Sections \ref{Sec:DefDeg}-\ref{Sec:Thm2} for its proof. As applications, in Section \ref{Sec3}, we outline how our degree theory can be used to prove the existence of solutions of the boundary Yamabe problem and the near-field reflector problem. In the Appendix, we collect some properties of the Laplace operator $\Delta : H^s \rightarrow H^{s-2}$ for $s \in [0,2]$ on a compact Riemannian manifold, which are needed in the body of the paper.

%++++++++++++++++%

\section{Statement of the main result}\label{Sec2}

In this section we introduce a degree theory for second order fully nonlinear elliptic operators with nonlinear oblique boundary conditions of general form,
	\begin{eqnarray}
		&& F[u] = f(\cdot,u,Du,D^2u),\quad \mbox{in }\Omega, \\
		&& G[u] = g(\cdot,u,Du),\quad \mbox{on }\partial\Omega,
	\end{eqnarray}
where $\Omega$ is a bounded smooth domain in Euclidean $n$-space, $\mathbb{R}^n$, and $f\in C^{3,\alpha}(\overline{\Omega}\times\mathbb{R}\times\mathbb{R}^n\times\mathcal{S}^n)$ and $g\in C^{4,\alpha}(\overline{\Omega}\times\mathbb{R}\times\mathbb{R}^n)$ are real valued functions, $0 < \alpha < 1$. Here $\mathcal{S}^n$ denotes the $n(n+1)/2$ dimensional linear space of $n\times n$ real symmetric matrices, and $Du=(D_iu)$ and $D^2u=[D_{ij}u]$ denote the gradient vector and Hessian matrix of the real valued function $u$.

Letting $(x,z,p,r)$ denote points in $\overline{\Omega}\times\mathbb{R}\times\mathbb{R}^n\times\mathcal{S}^n$, we shall adopt the following definitions of ellipticity and obliqueness for operators $F$ and $G$ \cite{GT}. An operator $F : C^{4,\alpha}(\overline{\Omega}) \to C^{2,\alpha}(\overline\Omega)$ is \emph{uniformly elliptic} on some bounded open subset $\mathcal{O}$ of $C^{4,\alpha}(\overline{\Omega})$ if there exists a constant $\lambda>0$ such that for all $u\in\mathcal{O}, x\in\overline{\Omega}$ and $\xi\in\mathbb{R}^n$ there holds
	\begin{equation}
		\frac{\partial f}{\partial r_{ij}}(x,u,Du,D^2u)\xi_i\xi_j \geq \lambda|\xi|^2.
	\end{equation}
An operator $G : C^{4,\alpha}(\overline{\Omega}) \to C^{3,\alpha}(\partial\Omega)$ is \emph{uniformly oblique} on $\mathcal{O}$ if there exists a constant $\chi>0$ such that for all $u\in\mathcal{O}$ and $x\in\partial\Omega$
	\begin{equation}
		\frac{\partial g}{\partial p}(x,u,Du)\cdot \gamma(x) \geq \chi,
	\end{equation}
where $\gamma(x)$ denotes the outer unit normal of $\partial\Omega$ at $x$. 

Let $\mathcal{O}\subset C^{4,\alpha}(\overline\Omega)$ be a bounded open set with $\partial\mathcal{O}\cap(F,G)^{-1}(0)=\emptyset$. Suppose that $F$ is uniformly elliptic on $\mathcal{O}$ and $G$ is uniformly oblique on $\mathcal{O}$. We will define an integer-valued degree for $(F,G)$ on $\mathcal{O}$ at $0$ along the line of \cite{Li,LL}.

\begin{theorem}\label{Thm:DegExistence}
There exists a unique integer-valued degree 
\[
deg: \Big\{((F,G),\mathcal{O}, 0): (F,G) \text{ and }\mathcal{O} \text{ are as above}\Big\} \rightarrow \ZZ
\]
which satisfies the following three properties:
\begin{enumerate}[(p1)]

\item $deg((F,G),\mathcal{O},0) = deg((F,G),\mathcal{O}_1,0) + deg((F,G),\mathcal{O}_2,0)$ whenever $\mathcal{O}_1$ and $\mathcal{O}_2$ are open disjoint subsets of $\mathcal{O}$ satisfying $(F,G)^{-1}(0) \cap (\mathcal{O}\setminus(\mathcal{O}_1 \cup \mathcal{O}_2)) = \emptyset$.

\item \emph{Homotopy invariance property}: If $t\mapsto (f_t,g_t)$ is continuous from $[0,1]$ to $C^{3,\alpha}(\overline\Omega\times\mathbb{R}\times\mathbb{R}^n\times\mathcal{S}^n)\times C^{4,\alpha}(\overline\Omega\times\mathbb{R}\times\mathbb{R}^n)$, $F_t$ is elliptic on $\mathcal{O}$, $G_t$ is oblique on $\mathcal{O}$ (both uniformly in $t \in [0,1]$), and $\partial\mathcal{O}\cap(F_t,G_t)^{-1}(0)=\emptyset$ for all $t\in[0,1]$, then $\deg((F_t,G_t),\mathcal{O},0)$ is independent of $t$. 

\item \emph{Compatibility with Leray-Schauder degree}: If $(F,G)$ is an invertible linear operator and $\mathcal{O}$ is a neighborhood of $0$, then
\[
\deg((F,G),\mathcal{O},0) = (-1)^{\dim E^{-}(F,G)}
\]
where
\begin{equation}
E^-(F,G) = \bigoplus_{\lambda_i < 0} \Big\{u \in C^{4,\alpha}(\bar\Omega): -(F[u],G[u]) = (\lambda_i u,0) \Big\}.
	\label{Eq:EminusDef}
\end{equation}

\end{enumerate}
\end{theorem}

As usual, the basic properties $(p1)$--$(p3)$ imply immediate consequences which we list below.

\begin{corollary}\label{Cor:a-d}
\begin{enumerate}[(a)]
\item If $\deg((F,G),\mathcal{O},0)\neq0$, then there exists $u\in\mathcal{O}$ such that $(F[u],G[u])=0$.

\item If $\overline{\mathcal{U}}\subset\mathcal{O}$ and $\overline{\mathcal{U}}\cap(F,G)^{-1}(0)=\emptyset$, then $\deg((F,G),\mathcal{O},0)=\deg((F,G),\mathcal{O}\setminus\overline{\mathcal{U}},0)$.

\item If $(F_1,G_1)|_{\partial\mathcal{O}} = (F_2,G_2)|_{\partial\mathcal{O}}$, then $\deg((F_1,G_1),\mathcal{O},0) = \deg((F_2,G_2),\mathcal{O},0)$. 

\item If $(F,G)(u_0)=0$, $(F,G)$ is Fr\'echet differentiable at $u_0$ and $(F',G')(u_0)$ is invertible. Then $\deg((F,G),\mathcal{O},0)=\deg((F',G'),\mathcal{B},0)$, where $\mathcal{O}$ is a neighborhood of $u_0$ in $C^{4,\alpha}(\overline\Omega)$ which does not contain any other points of $(F,G)^{-1}(0)$ and $\mathcal{B}$ is any bounded open set containing the origin.
\end{enumerate}
\end{corollary}

The next property is a mild extension of property $(p3)$.

\begin{corollary}\label{Cor:LSExt}
\begin{enumerate}[(a)]
\setcounter{enumi}{4}
\item Assume $(F,G)$ has linear leading terms, $(F,G)=(F_1,G_1)+(F_2,G_2)$ such that
	\begin{eqnarray*}
		&& (F_1,G_1)[u]=\left(a_{ij}(x)u_{ij}+b_i(x)u_i+c(x)u,\, (\beta_i(x)u_i+\ell(x)u)|_{\partial\Omega}\right) \\
		&& (F_2,G_2)[u]=\left(f_*(x,u,Du),\, b_*(x,u)|_{\partial\Omega}\right),
	\end{eqnarray*}
where $a_{ij}, b_i, c\in C^{3,\alpha}(\overline\Omega)$, $\beta_i, \ell\in C^{4,\alpha}(\partial\Omega)$, $f_*\in C^{3,\alpha}(\overline\Omega\times\mathbb{R}\times\mathbb{R}^n)$, and $b_*\in C^{4,\alpha}(\partial\Omega\times\mathbb{R})$. Assume that $F$ is elliptic, i.e. $(a_{ij})>0$ in $\overline\Omega$, and $G$ is oblique, i.e. $\beta\cdot\gamma>0$ on $\partial\Omega$, where $\gamma$ denotes the unit outer normal of $\partial\Omega$.

If $(F_1,G_1)$ is invertible, then for any open bounded $\mathcal{O}\in C^{4,\alpha}(\overline\Omega)$ such that $\partial\mathcal{O}\cap(F,G)^{-1}(0)=\emptyset$, there holds
	\begin{equation*}
		\deg((F,G),\mathcal{O},0)=(-1)^{\dim E^{-}(F_1,G_1)}\deg_{L.S.}(Id+(F_1,G_1)^{-1}\circ(F_2,G_2),\mathcal{O},0),
	\end{equation*}
where $E^{-}(F_1,G_1)$ is defined as in \eqref{Eq:EminusDef}.
\end{enumerate}
\end{corollary}

\section{Definition of the degree}\label{Sec:DefDeg}

Consider 
	\begin{eqnarray}
		S &:& C^{2,\alpha}(\overline\Omega) \to C^{\alpha}(\overline\Omega) \times C^{1,\alpha}(\partial\Omega)\nonumber \\
		  & & u \mapsto \left(\triangle u,\ \ (\gamma_iD_iu+u)|_{\partial\Omega}\right), 
	\label{Eq:SDef}
	\end{eqnarray}
where $\gamma$ is the outer unit normal of $\partial\Omega$, and 
	\begin{eqnarray}
		T &:& C^{3,\alpha}(\partial\Omega) \to C^{1,\alpha}(\partial\Omega) \nonumber\\
		  & & u \mapsto \triangle_T u-u, 
	\label{Eq:TDef}
	\end{eqnarray}
where $\triangle_T$ denotes the tangential Laplacian over $\partial\Omega$.	
It is well-known that $S, T$ are isomorphisms. 

Let $\tilde F, \tilde G$ be the composite maps as follows,
	\begin{eqnarray}
		&& \tilde F = \left(\tilde F_{(1)}, \tilde F_{(2)}\right) = S \circ F :  C^{4,\alpha}(\overline\Omega) \to C^{\alpha}(\overline\Omega) \times C^{1,\alpha}(\partial\Omega)\\
		&& \tilde G = T \circ G : C^{4,\alpha}(\overline\Omega) \to C^{1,\alpha}(\partial\Omega).
	\end{eqnarray}

Since $S$ and $T$ are isomorphisms, $(F,G)=0$ is equivalent to $(\tilde F,\tilde G)=0$. We are going to define a degree for $(F,G)$ by defining a degree for $(\tilde F,\tilde G)$.

As in \cite{Li}, we write
	\begin{eqnarray}
		&& \tilde F_{(1)}[u] = a_{st}(x,u,Du,D^2u)D_{iist}u + C_*(x,u,Du,D^2u,D^3u), \label{9}\\
		&& \tilde F_{(2)}[u] = \left.\left(a_{st}(x,u,Du,D^2u)D_{sti}u\gamma_i + E_*(x,u,Du,D^2u)\right)\right|_{\partial\Omega}, \label{10}\\
		&& \tilde G[u] = \left.\left(b_i(x,u,Du)\triangle_T(D_iu) + H_*(x,u,Du,D^2u)\right)\right|_{\partial\Omega},\label{11}
	\end{eqnarray}
where 
	\begin{eqnarray*}
		&& a_{st}(x,u,Du,D^2u) = \frac{\partial f}{\partial r_{st}}(x,u,Du,D^2u), \\
		&& b_i(x,u,Du) = \frac{\partial g}{\partial p_i}(x,u,Du).
	\end{eqnarray*}

To freeze coefficients in \eqref{9}--\eqref{11}, we make use of the following result, whose proof is postpone until Section \ref{Sec:Thm2}.

\begin{theorem}\label{Thm:LNBijection}
Let $a_{st}\in C^{1,\alpha}(\overline\Omega), b_i\in C^{1,\alpha}(\partial\Omega)$, where $1\leq i, s, t \leq n$. Assume that $(a_{st})$ is symmetric and there exists $\lambda>0$ such that $a_{st}(x)\xi_i\xi_j\geq\lambda|\xi|^2$ for all $\xi\in\mathbb{R}^n$ and $x\in\overline\Omega$, and there exists $\chi>0$ such that $b_i(x)\gamma_i(x)\geq\chi$ for all $x\in\partial\Omega$. For a constant $N$, define
	\begin{eqnarray*}
		L^N &:& C^{4,\alpha}(\overline\Omega) \to C^\alpha(\overline\Omega) \times C^{1,\alpha}(\partial\Omega) \times C^{1,\alpha}(\partial\Omega) \\
		 && w \mapsto \left(L^N_{(1)}w,\ L^N_{(2)}w,\ L^N_{(3)}w\right), \nonumber
	\end{eqnarray*}
where
	\begin{eqnarray*}
		&& L^N_{(1)}w = a_{st}D_{iist}w - N a_{st}D_{st}w, \\
		&& L^N_{(2)}w = \left.\left(a_{st}D_{sti}w\gamma_i\right)\right|_{\partial\Omega}, \\
		&& L^N_{(3)}w = \left.\left(b_i\triangle_T(D_iw) - N b_iD_iw-Nw\right)\right|_{\partial\Omega}.
 	\end{eqnarray*}
	
Then there exists some constant $N_0$, depending only on $\|a_{st}\|_{C^{1,\alpha}}, \|b_i\|_{C^{1,\alpha}}, n,\lambda,\chi$ such that $L^N$ is an isomorphism for all $N>N_0$. Furthermore, $L^N$ depends continuously on $a_{st}, b_i$ with respect to the corresponding topologies. 
\end{theorem}

We are in a position to define an integer-valued degree for $(F,G) : \mathcal{O}\to C^{2,\alpha}(\overline\Omega)\times C^{3,\alpha}(\partial\Omega)$. Note that $\tilde F, \tilde G$ in (9)--(11) can be represented as 
	\begin{equation}
		(\tilde F[u], \tilde G[u]) = L^{u,N}[u] + R^{u,N}[u], 
	\label{Eq:LRDef}
	\end{equation}
where $L^{u,N}=(L^{u,N}_{(1)},L^{u,N}_{(2)},L^{u,N}_{(3)})$ and  $R^{u,N}=(R^{u,N}_{(1)},R^{u,N}_{(2)},R^{u,N}_{(3)})$ are
	\begin{eqnarray*}
		&& L^{u,N}_{(1)}w=a_{st}(x,u,Du,D^2u)D_{iist}w-Na_{st}(x,u,Du,D^2u)D_{st}w, \\
		&& L^{u,N}_{(2)}w=\left.\left(a_{st}(x,u,Du,D^2u)D_{sti}w\gamma_i\right)\right|_{\partial\Omega}, \\
		&& L^{u,N}_{(3)}w=\left.\left(b_i(x,u,Du)\triangle_T(D_iw)-Nb_i(x,u,Du)D_iw-Nw\right)\right|_{\partial\Omega}, \\
		&& R^{u,N}_{(1)}w=Na_{st}(x,u,Du,D^2u)D_{st}w+C_*(x,u,Du,D^2u,D^3u),\\
		&& R^{u,N}_{(2)}w=E_*(x,u,Du,D^2u)|_{\partial\Omega},\\
		&& R^{u,N}_{(3)}w=Nb_i(x,u,Du)D_iw+Nw+H_*(x,u,Du,D^2u)|_{\partial\Omega}.
	\end{eqnarray*}

One can see that $R^{u,N}[u]$ maps $C^{4,\alpha}(\overline\Omega)$ into $C^{1,\alpha}(\overline\Omega)\times C^{2,\alpha}(\partial\Omega)\times C^{2,\alpha}(\partial\Omega)$.
According to Theorem \ref{Thm:LNBijection}, there exists some positive number $N_0$ such that $L^{u,N}$ is an isomorphism for any $N>N_0$.
By \cite[Theorem 7.3]{ADN}, $(L^{u,N})^{-1}$ maps $C^{1,\alpha}(\overline\Omega)\times C^{2,\alpha}(\partial\Omega)\times C^{2,\alpha}(\partial\Omega)$ into $C^{5,\alpha}(\overline\Omega)$, and its norm as a linear map between these spaces is bounded by a constant depends only on $\|a_{st}\|_{C^{2,\alpha}}, \|b_i\|_{C^{3,\alpha}}, \lambda$ and $\chi$. It follows that
	\begin{equation*}
		u \mapsto (L^{u,N})^{-1}R^{u,N}[u]
	\end{equation*}
is a compact operator from $\mathcal{O}$ to $C^{4,\alpha}(\overline\Omega)$. 

Moreover, $(\tilde F, \tilde G)[u]=0$ is the same as $u+(L^{u,N})^{-1}R^{u,N}[u]=0$, i.e.,
	\begin{equation*}
		\partial\mathcal{O}\cap(Id+(L^{u,N})^{-1}R^{u,N})^{-1}(0)=\partial\mathcal{O}\cap(F,G)^{-1}(0)=\emptyset.
	\end{equation*}
Therefore, we can define the degree of $(F,G)$ as the Leray-Schauder degree of the map $u \mapsto u+(L^{u,N})^{-1}R^{u,N}[u]$. (See e.g. \cite{NirBook} for the definition of the Leray-Schauder degree.) More precisely we have the following definition.

\begin{definition}\label{Def:DegDef}
Let $(F,G)$ and $\mathcal{O}$ be as in Theorem \ref{Thm:DegExistence}. We define a degree of $(F,G)$ on $\mathcal{O}$ at $0$ by
	\begin{equation}\label{new def}
		\deg\left((F,G),\mathcal{O},0\right)=\deg_{L.S.}(Id+(L^{u,N})^{-1}R^{u,N},\mathcal{O},0),
	\end{equation}
where $N>N_0$ and $N_0$ is the constant in Theorem \ref{Thm:LNBijection}.
\end{definition}

Note that, by the homotopy invariance of the Leray-Schauder degree, the degree above is independent of $N>N_0$. In Sections \ref{Sec:BCE} and \ref{Sec:Thm1} below, we shall prove Theorem \ref{Thm:LNBijection} and verify the above defined degree \eqref{new def} satisfies the required properties $(p1)$--$(p3)$ in Theorem \ref{Thm:DegExistence}.

\section{Some boundary estimates}\label{Sec:BCE}

We start with some boundary estimates for linear elliptic systems with oblique boundary conditions. We use $B_r^+$ to denote $\{x \in \RR^n: |x| < r, x_n > 0\}$ and let $\Gamma = \{x_n = 0\} \cap B_1$. Fix some integers $m \geq 1$ and $n \geq 2$. In the sequel, repeated Roman indices are summed from $1$ to $n$ and repeated Greek indices are summed from $1$ to $m$. Consider the system 
	\begin{eqnarray}
		&& a_{st}^{\alpha\beta} D_{st}u^\beta + d_s^{\alpha\beta}\,D_s u^\beta + c^{\alpha\beta}\,u^\beta = f^\alpha \quad\mbox{in }B_1^+, \label{B21}\\
		&& b_i^{\alpha\beta} D_i u^\beta = g^\alpha \quad\mbox{on } \Gamma.\label{B22}
	\end{eqnarray}
We assume that $a_{st}^{\alpha\beta}$ is uniformly strongly elliptic in $B_1^+$, i.e. there exists $\lambda > 0$ such that
\[
a_{st}^{\alpha\beta}(x)\,\xi_s^\alpha\,\xi_t^\beta \geq \lambda|\xi|^2 \text{ for all } x \in B_1^+ \text{ and } \xi \in \RR^{n \times m},
\]
and $b_i^{\alpha\beta}$ is uniformly oblique along $\Gamma$, i.e. there exists $\chi > 0$ such that
\[
-b_n^{\alpha\beta}(x)\eta^\alpha\,\eta^\beta \geq \chi|\eta|^2 \text{ for all } x \in \Gamma \text{ and } \eta \in \RR^m.
\]

\begin{lemma}\label{Lem:LocalBdryEst}
For $m \geq 1$, assume that $a_{st}^{\alpha\beta} \in W^{1,\infty}(\bar B_1^+)$ is uniformly elliptic and symmetric, $d_s^{\alpha\beta} \in L^\infty(B_1^+)$, $c^{\alpha\beta} \in L^\infty(B_1^+)$ and $b_i^{\alpha\beta} \in W^{1,\infty}(\Gamma)$ is uniformly oblique along $\Gamma$. Let $F = (f^\alpha)_{\alpha=1}^m \in L^2(B_1^+)$, $G = (g^\alpha)_{\alpha=1}^m \in L^2(\Gamma)$, and $U = (u^\alpha)_{\alpha = 1}^m \in H^2(B_1^+)$ be a solution of the oblique boundary value problem \eqref{B21}-\eqref{B22}. Then we have the estimate
	\begin{equation}
		\|DU\|_{L^2(B_{1/2}^+)} + \|DU\|_{L^2(\Gamma \cap B_{1/2})} \leq C \Big[\|F\|_{L^2(B_1^+)} + \|G\|_{L^2(\Gamma)}  + \|U\|_{L^2(B_1^+)}\Big],
	\end{equation}
where $C$ depends only on $\|a_{st}^{\alpha\beta}\|_{W^{1,\infty}(\bar B_1^+)}$, $\|d_s^{\alpha\beta}\|_{L^\infty(B_1^+)}$, $\|c^{\alpha\beta}\|_{L^\infty(B_1^+)}$, $\|b_i^{\alpha\beta}\|_{W^{1,\infty}(\Gamma)}$ and the ellipticity and obliqueness constants $\lambda$ and $\chi$.
\end{lemma}

\begin{proof} We adapt the proof of the well-known Rellich identity for harmonic functions. 

Fix $0 < r_1 < r_2 < 1$. In the sequel, $C$ will denote some positive constant that may very from line to line and depends only on the coefficients of the equation.

Select a smooth cutoff function $\varphi$ satisfying
	\begin{equation*}
	\left\{\begin{array}{ll}
		\varphi\equiv1 \mbox{ in }B_{r_1}(0), &\\ \varphi\equiv0 \mbox{ on }\mathbb{R}^n-B_{(r_1 + r_2)/2}(0), \\
		0\leq\varphi\leq1, |D^k \varphi| \leq C\,A^k, &
	\end{array}\right.
	\end{equation*}
where
\[
A = \frac{1}{r_2 - r_1} > 1.
\]

\noindent 1. Testing \eqref{B21} against $\varphi^2U $, we have
	\begin{equation}
		\int_{B_{r_1}^+}|DU|^2 \leq\ C\int_{B_{r_2}^+} |F|^2 + C\,A^2\,\int_{B_{r_2}^+}|U|^2 + C\,K\int_{\Gamma \cap B_{r_2}}|U|^2 + K^{-1} \int_{\Gamma \cap B_{r_2}}|DU|^2 
			\label{B27}
	\end{equation}
for large $K >0$.

\medskip
\noindent 2. Testing \eqref{B21} against $\varphi^2 D_n U$, we have
	\begin{equation*}
	\begin{split}
		-\int_{B_1^+} \varphi^2\,f^\alpha\,D_n u^\alpha 
		  \leq& \int_{\Gamma}\varphi^2 D_n u^\alpha D_s u^\beta a_{sn}^{\alpha\beta}
			+ \int_{B_1^+} \varphi^2 a_{st}^{\alpha\beta}  D_s u^\beta D_{nt} u^\alpha \\ 
			&\qquad
				+ C\,A\,\int_{B_1^+}\varphi |DU|^2 + C\,\int_{B_1^+} \varphi^2\,|U|^2\\
		 \leq& \int_{\Gamma}\varphi^2 D_n u^\alpha D_s u^\beta a_{sn}^{\alpha\beta}
			+ \int_{B_1^+}\varphi^2 D_n(\frac12 a_{st}^{\alpha\beta}  D_s u^\beta D_t u^\alpha) 
			\\
			&\qquad 
				+ C\,A\,\int_{B_1^+}\varphi |DU|^2  + C\,\int_{B_1^+} \varphi^2\,|U|^2\\
		 \leq& \int_{\Gamma}\varphi^2 D_n u^\alpha D_s u^\beta a_{sn}^{\alpha\beta}
			- \frac1C\int_{\Gamma} \varphi^2|DU|^2 + C\,A\,\int_{B_1^+}\varphi |DU|^2 +   C\,\int_{B_1^+} \varphi^2\,|U|^2,
	\end{split}	
	\end{equation*}
thus, 
	\begin{equation}
	\begin{split}
		\int_{\Gamma \cap B_{r_1}} |DU|^2 
			\leq &C \max\Big\{0, \int_{\Gamma}\varphi^2 D_n u^\alpha D_s u^\beta a_{sn}^{\alpha\beta} \Big\} + C\int_{B_{r_2}^+} |F|^2 + C\,A\int_{B_{r_2}^+}|DU|^2  + C\,\int_{B_{r_2}^+} \,|U|^2.	
	\end{split}
	\label{B29}	
	\end{equation}

\noindent 3. Let $\tilde b_n^{\beta\alpha}$ be the inverse of $b_n^{\alpha\beta}$, i.e. $\tilde b_n^{\beta\alpha}\,b_n^{\alpha\gamma} = \delta^{\beta\gamma}$. Set $v^\alpha  = \tilde b_n^{\alpha\gamma}\sum_{i = 1}^{n-1} b_i^{\gamma\mu}\,D_iu^\mu$. Then
\[
D_n u^\alpha = \tilde b_n^{\alpha\beta}\,g^\beta - v^\alpha. 
\]
Testing \eqref{B21} against $\varphi^2 V = (\varphi^2 v^\alpha)_{\alpha = 1}^m$, we have
	\begin{equation*}
	\begin{split}
		-\int_{B_1^+} \varphi^2\,f^\alpha\,v^\alpha 
		=& -\int_{B_1^+}(\varphi^2 v^\alpha a_{st}^{\alpha\beta})D_{st} u^\beta
			-\int_{B_1^+}\varphi^2 v^\alpha\,d_k^{\alpha\beta} D_k u^\beta \\
		  =& \int_{\Gamma}\varphi^2 D_s u^\beta v^\alpha a_{sn}^{\alpha\beta} + \int_{B_1^+} \varphi^2 a_{st}^{\alpha\beta}  D_t v^\alpha  D_s u^\beta\\
			&\qquad + O\Big(A\int_{B_{r_2}^+} |DU|^2\,dx + \int_{B_{r_2}^+} |U|^2\,dx\Big)
	\end{split}
	\end{equation*}
Noting that
\[
a_{st}^{\alpha\beta}  D_t v^\alpha  D_s u^\beta = \sum_{i = 1}^{n-1} a_{st}^{\alpha\beta}\,D_t (\tilde b_n^{\alpha\gamma}\,b_i^{\gamma\,\mu}\,D_i u^\mu)\,D_s u^\beta
\]

This implies that
	\begin{equation*}
		 \left|\int_{\Gamma}\varphi^2 D_s u^\alpha v^\alpha a_{sn}\right|  \leq CA^{-1}\int_{B_{r_2}^+} |F|^2 + C\,A\, \int_{B_{(r_1 + r_2)/2}^+}|DU|^2.
	\end{equation*}
Hence, by the uniform obliqueness we obtain in case the integral on the left hand side below is positive that 
	\begin{equation*}
	\begin{split}
		\int_{\Gamma}\varphi^2 D_s u^\alpha D_n u^\alpha a_{sn} &\leq C\int_{\Gamma}\varphi^2D_s u^\alpha (-b_n) D_n u^\alpha a_{sn} = -C \int_{\Gamma}\varphi^2 D_s u^\alpha \big(g^\alpha -v^\alpha\big) a_{sn}\\
				&\leq \varepsilon\int_{\Gamma} \varphi^2 |DU|^2 
					+ C_\varepsilon\int_{\Gamma \cap B_{r_2}} |G|^2\\
				&\qquad +  CA^{-1}\int_{B_{r_2}^+} |F|^2 + C\,A\, \int_{B_{(r_1 + r_2)/2}^+}|DU|^2
	\end{split}
	\end{equation*}
for any $\varepsilon>0$ small. Recalling \eqref{B29} we obtain
	\begin{equation}
	\begin{split}
		\int_{\Gamma \cap B_{r_1}} |DU|^2 
			\leq & CA^{-1}\int_{B_{r_2}^+} |F|^2 + C\int_{\Gamma \cap B_{r_2}} |G|^2 + C\,A\int_{B_{(r_1 + r_2)/2}^+}|DU|^2.	
	\end{split}
	\label{B29x}	
	\end{equation}
Recalling \eqref{B27} with $K = \varepsilon^{-1}\,A$ for some small $\varepsilon > 0$, we obtain
	\begin{equation}
	\begin{split}
		\int_{\Gamma \cap B_{r_1}} |DU|^2 
			\leq & CA^{-1}\int_{B_{r_2}^+} |F|^2 + C\int_{\Gamma \cap B_{r_2}} |G|^2  \\
			&\qquad
+ C\,A^3\,\int_{B_{r_2}^+}|U|^2 + C\,\epsilon^{-1}\,A^2\,\int_{\Gamma \cap B_{r_2}}|U|^2 + \varepsilon\,\int_{\Gamma \cap B_{r_2}}|DU|^2.
	\end{split}
	\label{B27y}
	\end{equation}

\noindent 4. Let 
\[
\Phi(r) = \int_{\Gamma \cap B_{r}} |DU|^2.
\]
Then by \eqref{B27y} with $\varepsilon = \frac{1}{2}$,
	\begin{equation}
	\begin{split}
		\Phi(r_1)
			\leq & \frac{1}{2}\Phi(r_2) + C\int_{B_1^+} |F|^2 + C\int_{\Gamma} |G|^2  
+ C\,A^3\,\int_{B_1^+}|U|^2 + C\,A^2\,\int_{\Gamma}|U|^2\\
			\leq & \frac{1}{2}\Phi(r_2) + \frac{C}{(r_2 - r_1)^3}\Big[\int_{B_1^+} |F|^2 + C\int_{\Gamma} |G|^2  
+ \int_{B_1^+}|U|^2 + \int_{\Gamma}|U|^2\Big]\\
			=: &\frac{1}{2}\Phi(r_2) + \frac{C}{(r_2 - r_1)^3}X.
	\end{split}
	\label{B27z}
	\end{equation}
A standard iteration (see e.g. \cite[Lemma 1.1]{GiaGiusti}) leads to
\[
\Phi(\frac{1}{2}) \leq CX.
\]
In other words,

\[
\int_{\Gamma \cap B_{1/2}} |DU|^2 \leq C\int_{B_1^+} |F|^2 + C\int_{\Gamma} |G|^2  
+ C\,\int_{B_1^+}|U|^2 + C\,\int_{\Gamma}|U|^2.
\]
Combining with \eqref{B27} we conclude the proof.
\end{proof}

As a consequence, we obtain the following boundary estimate for scalar oblique boundary value problems.

\begin{lemma}\label{Lem:Lem1}
Assume that $a_{st} \in W^{1,\infty}(\Omega)$ is uniformly elliptic and symmetric and $b_i \in W^{1,\infty}(\partial\Omega)$ is uniformly oblique along $\partial\Omega$. Let $g \in L^2(\partial\Omega)$ and let $w\in H^3(\Omega)$ be a solution of the oblique boundary value problem 
	\begin{eqnarray}
		&& a_{st}D_{st}w = 0\quad\mbox{in }\Omega, \label{21}\\
		&& b_iD_iw+w = g\quad\mbox{on }\partial\Omega.\label{22}
	\end{eqnarray}
Then we have the estimates
	\begin{eqnarray}
		\|w\|_{L^2(\partial\Omega)} + \|Dw\|_{L^2(\partial\Omega)} &\leq& C \|g\|_{L^2(\partial\Omega)}, \label{Lemma1est1}\\
	        \|D^2 w\|_{L^2(\partial\Omega)} &\leq& C\|g\|_{H^1(\partial\Omega)},
	         \label{Lemma1est2X}
	\end{eqnarray}
where the constant $C$ depends only on $\|a_{st}\|_{W^{1,\infty}(\Omega)}$, $\|b_i\|_{W^{1,\infty}(\partial\Omega)}$ and the ellipticity and obliqueness constants $\lambda$ and $\chi$.
\end{lemma}

\begin{remark}
Later on, we will use the following consequence of \eqref{Lemma1est2X}:
\begin{eqnarray}
	        \|\nabla_T w_i\|_{L^2(\partial\Omega)} &\leq&
	        	        C\Big( \|\sum_j b_j \nabla_T w_j\|_{L^2(\partial\Omega)} +  \|g\|_{L^2(\partial\Omega)} \Big)
	         \label{Lemma1est2}
	\end{eqnarray}
for all $i = 1, \ldots, n$, where $\nabla_T$ denotes the covariant derivative along $\partial\Omega$.
\end{remark}

\begin{proof}
The proof is based on the local boundary estimate in Lemma \ref{Lem:LocalBdryEst}.

For every point $x_0 \in \partial \Omega$, we can find some sufficiently small $r > 0$ and a diffeomorphism $\Phi: \Omega\cap B_r(x_0) \rightarrow B_1^+(0)$.

Define
	\begin{equation}
		\tilde w(y) := w\left(\Phi^{-1}(y)\right),\quad (y\in B_1^+).	
	\label{34}
	\end{equation}
It is straightforward to check that $\tilde w$ satisfies
	\begin{equation}\label{neweqn}
		\left\{\begin{array}{l} \tilde a_{st}D_{st} \tilde w + \tilde d_s\,D_s \tilde w = 0 \text{ in } B_1^+,\\
		 \tilde b_i \,D_i \tilde w + \tilde w = \tilde g \text{ on } \Gamma = \{y_n = 0\} \cap B_1.
		 \end{array}\right.
	\end{equation}
where $\tilde a = (D\Phi)^T\,a\,D\Phi$ is uniformly elliptic, $\tilde b = (D\Phi)^T b$ is uniformly oblique, and $\tilde d = a_{st}\partial_{st} \Phi$.

\medskip
\noindent\underline{Proof of \eqref{Lemma1est1}:} Apply Lemma \ref{Lem:LocalBdryEst} to \eqref{neweqn} by setting $m=1, F=0$ and $G=g-\tilde w$, we have
	\begin{equation*}
		\int_{\Gamma\cap B_{1/2}}|D\tilde w|^2 + \int_{B_{1/2}^+} |D\tilde w|^2 \leq C\left(\int_\Gamma|g|^2+\int_\Gamma|\tilde w|^2+\int_{B_1^+}|\tilde w|^2\right).
	\end{equation*}
Consequently
	\begin{equation}
		\int_{\partial\Omega\cap B_{r/2}}|D w|^2 + \int_{\Omega\cap B_{r/2}}|D w|^2 \leq  C\left(\int_{\partial\Omega\cap B_r}|g|^2+\int_{\partial\Omega\cap B_r}|w|^2+\int_{\Omega\cap B_r}|w|^2\right).
	\end{equation}

Since $\partial\Omega$ is compact, we can cover $\partial\Omega$ with finitely many small balls $B_r(x_i)$, $x_i\in\partial\Omega$, $i=1,\cdots,N$. Summing the estimates for $i$ from $1$ to $N$, we obtain that
	\begin{equation*}
		\int_{\partial\Omega}|Dw|^2 + \int_{\Omega_r} |Dw|^2 \leq C\left(\int_{\partial\Omega}|g|^2+\int_{\partial\Omega}|w|^2+\int_{\Omega_r}|w|^2\right),
	\end{equation*}
for a small constant $r>0$, where $\Omega_r=\{x\in\Omega\,:\,\dist(x,\partial\Omega)<r\}$. Together with interior $L^2$ gradient estimates for \eqref{21}, we arrive at
	\begin{equation}
		\int_{\partial\Omega}|Dw|^2 + \int_{\Omega} |Dw|^2 \leq C\left(\int_{\partial\Omega}|g|^2+\int_{\partial\Omega}|w|^2+\int_{\Omega}|w|^2\right).
		\label{45}
	\end{equation}
We can then apply a standard argument using compactness and the uniqueness of the problem \eqref{21}--\eqref{22} to simplify estimate \eqref{45} to
	\begin{equation}
		\int_{\partial\Omega}|Dw|^2 + \int_{\Omega} |Dw|^2 \leq C\int_{\partial\Omega}|g|^2.
	\label{21I13.01}
	\end{equation}
This finishes the proof of \eqref{Lemma1est1}.

\medskip
\noindent\underline{Proof of \eqref{Lemma1est2X}:} As before, we first investigate \eqref{neweqn}. Fix some $\tau = 1, \ldots, n - 1$. By differentiating equation \eqref{neweqn} with respect to $x_\tau$, one has
	\begin{equation}\label{neweqndif}
		\tilde a_{st}D_{st}(D_\tau \tilde w) + \tilde d_kD_k(D_\tau \tilde w) = -D_\tau(\tilde a_{st})D_{st} \tilde w-(D_\tau \tilde d_k)D _k \tilde w.
	\end{equation}
Also, from \eqref{neweqn}, we can write $D_{nn}\tilde w$ as a combination of $\{D_{st} \tilde w: (s,t) \neq (n,n)\}$ and $\{D_k \tilde w: 1 \leq k \leq n\}$. Thus, $W=(D_\tau \tilde w)_{\tau=1}^{n-1}$ satisfies
\begin{equation}\label{nedMod}
\tilde a_{st}D_{st} W^\tau + \tilde d_k^{\tau\tau'} D_k W^{\tau'} = f^\tau 
\end{equation}
for some smooth $d_k^{\tau\tau'}$ and $F = (f^\tau)_{\tau = 1}^{n-1}$ satisfying $|F| \leq C|D\tilde w|$.

Applying Lemma \ref{Lem:LocalBdryEst} to \eqref{nedMod}, we have
\[
   \|DW\|_{L^2(\Gamma\cap B_{1/2})} \leq C\left[\|F\|_{L^2(B_1^+)}+\|b_i\,D_i W\|_{L^2(\Gamma)} +\|W\|_{L^2(B_1^+)}+\|W\|_{L^2(\Gamma)}\right].
 \]
Returning to $w$ and using the compactness of $\partial\Omega$ and estimate \eqref{21I13.01}, we obtain that
\[
	  \int_{\partial\Omega}|D\nabla_T w|^2 \leq C \left(\int_{\partial\Omega}|\sum_j b_j \nabla_T w_j|^2 + \int_{\partial\Omega}|g|^2\right).
\]
By the equation, $D_{\gamma\gamma} w$ is also under control. \eqref{Lemma1est2X} is proved.
\end{proof}

%++++++++++++++++%
\section{Proof of Theorem \ref{Thm:DegExistence}}\label{Sec:Thm1}

We assume for now the correctness of Theorem \ref{Thm:LNBijection}, whose proof will be carried out in the next section. Then the degree $\deg((F,G),\mathcal{O}, 0)$ in Definition \ref{Def:DegDef} is well-defined. Properties $(p1)$--$(p2)$ follow from the properties of the Leray-Schauder degree. For $(p3)$, we prove the more general statement in Corollary \ref{Cor:LSExt}. To this end we use the following lemma on the semi-finiteness of a linear operator.

\begin{lemma}\label{SemiFiniteness}
Assume $a_{ij}$ $\in$ $W^{1,\infty}(\Omega)$, $b_i$, $c$ $\in$ $L^\infty(\Omega)$, $\beta_i$ $\in$ $W^{1,\infty}(\partial\Omega)$ and $\ell$ $\in$ $L^\infty(\partial\Omega)$. Assume furthermore that $(a_{ij})$ $>$ $\lambda\,I$ in $\bar\Omega$ and $\beta \cdot \gamma$ $>$ $\chi$ on $\partial\Omega$ for some positive constants $\lambda$ and $\chi$. Then there exists $\mu_*$ depending on  $\|a_{ij}\|_{W^{1,\infty}(\Omega)}$, $\|b_i\|_{L^\infty(\Omega)}$, $\|c^+\|_{L^\infty(\Omega)}$, $\|\beta_i\|_{W^{1,\infty}(\partial\Omega)}$, $\|\ell^-\|_{L^\infty(\partial\Omega)}$, $\lambda$ and $\chi$ such that for any $\mu$ $>$ $\mu_*$, the problem 
\begin{equation}
\left\{\begin{array}{l}
a_{ij}(x)\,u_{ij} + b_i(x)\,u_i + c(x)\,u = \mu u \text{ in } \Omega,\\
\beta_i(x)\,u_i + \ell(x)\,u = 0 \text{ on } \partial\Omega
\end{array}\right.
\label{LinearEVP}
\end{equation}
has no non-trivial solution in $H^2(\Omega)$.
\end{lemma}

\begin{remark}
If $\ell$ $\geq$ $0$ on $\partial\Omega$ and the coefficients are smooth, the result follows directly from the maximum principle. In fact, $\mu_*$ can then be taken to be $\|c^+\|_{C^0(\bar\Omega)}$.
\end{remark}

\begin{proof} We use energy method. Assume that $u$ $\in$ $H^1(\Omega)$ is a solution to \eqref{LinearEVP}. We will use $C$ to denote some positive constant which may vary from lines to lines and depends only on $\|a_{ij}\|_{W^{1,\infty}(\Omega)}$, $\|b_i\|_{L^\infty(\Omega)}$, $\|c^+\|_{L^\infty(\Omega)}$, $\|\beta_i\|_{W^{1,\infty}(\partial\Omega)}$, $\|\ell^-\|_{L^\infty(\partial\Omega)}$, $\lambda$ and $\chi$. In particular, $C$ is always independent of $\mu$.

Multiplying the first equation in \eqref{LinearEVP} by $u$ then integrating over $\Omega$, we get
\begin{align*}
\mu\int_\Omega u^2\,dx
	&\leq -C^{-1}\int_\Omega |D u|^2\,dx + C\int_\Omega u^2\,dx + \int_{\partial\Omega} u\,a_{ij}\,u_i\,\gamma_j\,d\sigma(x).
\end{align*}
To proceed, we write 
\[
a_{ij}\gamma_j = p \beta_i  + X_i
\]
where $p = \frac{a_{ij}\gamma_i\gamma_j}{\beta \cdot \gamma}$ and $X\cdot \gamma = 0$. Note that $p > 0$ is bounded thanks to ellipticity and obliqueness. It follows that
\begin{align}
C^{-1}\int_\Omega |D u|^2\,dx + (\mu - C)\int_\Omega u^2\,dx
	&\leq \int_{\partial\Omega} u\,\Big[p\,\beta_i u_i + X_i\,u_i\Big]\,d\sigma(x)\nonumber\\
	&= \int_{\partial \Omega} u\,[-p\,\ell\,u + X(u)]\,d\sigma(x).\label{LSF::Bound1}
\end{align}

Note that, by Stoke's theorem, we have
\[
\int_{\partial\Omega} u\,X(u)\,d\sigma(x) = \frac{1}{2}\int_{\partial\Omega} X(u^2)\,d\sigma(x) \leq C \int_{\partial\Omega} u^2\,d\sigma(x).
\]
Returning to \eqref{LSF::Bound1} we hence get
\begin{align*}
C^{-1}\int_\Omega |D u|^2\,dx + (\mu - C)\int_\Omega u^2\,dx
	&\leq C\int_{\partial \Omega} u^2\,d\sigma(x)\\
	&\leq \epsilon \int_\Omega |D u|^2\,dx + C_\epsilon\int_\Omega u^2\,dx
\end{align*}
for any small $\epsilon > 0$. Here we have used the compactness of the embedding $H^1(\Omega)$ $\hookrightarrow$ $L^2(\partial\Omega)$. The assertion follows by choosing $\epsilon = \frac{1}{2C}$.
\end{proof}

\begin{proof}[Proof of Corollary \ref{Cor:LSExt}]
As before, set $(\tilde F,\tilde G)$ $=$ $(S\circ F,T \circ G):$ $C^{4,\alpha}(\bar\Omega)$ $\rightarrow$ $C^\alpha(\bar\Omega) \times C^{1,\alpha}(\partial\Omega) \times C^{3,\alpha}(\partial\Omega)$, where $S$ and $T$ are given by \eqref{Eq:SDef} and \eqref{Eq:TDef}. 

In case $(F,G)$ has the above special form, the operator $L = L_{u,N}: C^{4,\alpha}(\bar\Omega) \rightarrow C^\alpha(\bar \Omega) \times C^{1,\alpha}(\partial\Omega) \times C^{1,\alpha}(\partial\Omega)$ defined in \eqref{Eq:LRDef} takes the form
\[
Lw = \Big(a_{st} D_{iist} w - Na_{st}D_{st} w, (a_{st} D_{ist} w\,\gamma_i)\Big|_{\partial\Omega}, (b_i\,\Delta_T(D_i w) - N\,b_i\,D_i w - N w)\Big|_{\partial\Omega}\Big).
\]
By Theorem \ref{Thm:LNBijection}, we can select $N$ sufficiently large such that $L$ is invertible, $L^{-1} \circ (\tilde F, \tilde G): C^{4,\alpha}(\bar\Omega) \rightarrow C^{4,\alpha}(\bar\Omega)$ is of the form $\textrm{Id} + \textrm{Compact}$ and
\[
\deg((F,G),\mathcal{O},0) = \deg_{L.S.}(L^{-1}\circ (\tilde F,\tilde G),\mathcal{O},0).
\]

Set $(\tilde F_1, \tilde G_1) = (S\circ F_1, T \circ G_1): C^{4,\alpha}(\bar\Omega) \rightarrow C^\alpha(\bar\Omega) \times C^{1,\alpha}(\partial\Omega) \times C^{3,\alpha}(\partial\Omega)$. By our hypotheses, $(\tilde F_1, \tilde G_1)$ is invertible. Thus, by the product rule of the Leray-Schauder degree, 
\[
\deg((F,G),\mathcal{O},0) = \sum_{\mathcal{U}}\deg_{L.S.}(L^{-1}\circ (\tilde F_1,\tilde G_1),\mathcal{U},0)\,\deg_{L.S.}((\tilde F_1, \tilde G_1)^{-1}\circ (\tilde F,\tilde G),\mathcal{O},\mathcal{U}),
\]
where the summation is made over the connected components of $C^{4,\alpha}(\bar\Omega) \setminus (F_1,G_1)^{-1}\circ(F,G)(\partial\mathcal{O})$. It is evident that $\deg_{L.S.}(L^{-1}\circ (\tilde F_1, \tilde G_1),\mathcal{U},0)$ $=$ $0$ if $0$ $\notin$ $\mathcal{U}$. Hence
\[
\deg((F,G),\mathcal{O},0) = \deg_{L.S.}(L^{-1}\circ (\tilde F_1, \tilde G_1),\mathcal{\tilde O},0)\,\deg_{L.S.}((\tilde F_1, \tilde G_1)^{-1}\circ (\tilde F,\tilde G),\mathcal{O},0),
\]
where $\mathcal{\tilde O}$ is the connected component of $C^{4,\alpha}(\bar\Omega) \setminus (F_1,G_1)^{-1}\circ(F,G)(\partial\mathcal{O})$ containing $0$. As $(\tilde F_1, \tilde G_1)^{-1}\circ (\tilde F,\tilde G)$ $=$ $( F_1, G_1)^{-1}\circ (F, G)$ $=$ $Id + (F_1,G_1)^{-1}\circ(F_2,G_2)$, it remains to show that
\begin{equation}
d := \deg_{L.S.}(L^{-1}\circ (\tilde F_1,\tilde G_1),\mathcal{\tilde O},0) = (-1)^{\dim E^-(F_1,G_1)}.
\label{LL::Eqn01}
\end{equation}

Define $(F_3, G_3):$ $C^{4,\alpha}(\bar\Omega)$ $\rightarrow$ $C^{2,\alpha}(\bar\Omega) \times C^{1,\alpha}(\partial\Omega)$ by
\[
(F_3[w], G_3[w]) = \Big(a_{ij}\,w_{ij},(\beta_i\,w_i + w)\Big|_{\partial\Omega}\Big).
\]
For $0$ $\leq$ $t$ $\leq$ $1$, define $L_t:$ $C^{4,\alpha}(\bar\Omega)$ $\rightarrow$ $C^\alpha(\bar\Omega) \times C^{1,\alpha}(\partial\Omega) \times C^{1,\alpha}(\partial\Omega)$ by
\[
\begin{split}
L_t\,w 	
	&= \Big((1-t)\,a_{st}\,w_{iist} + t\Delta F_3[w] - N\,F_3[w], \Big((1-t)a_{st}\,w_{sti}\,\gamma_i + t\,\frac{\partial F_3[w]}{\partial \gamma}\Big)\Big|_{\partial\Omega},\\
		&\qquad\qquad \Big((1-t) b_i \Delta_T (D_i w) + t\Delta_T G_3[w] - NG_3[w]\Big)\Big|_{\partial\Omega}\Big).
\end{split}
\]
We can then apply the proof of Theorem \ref{Thm:LNBijection} to see that, for sufficiently large $N$, $L_t$ is an isomorphism for each $t$ $\in$ $[0,1]$. Furthermore, as $L_t - (\tilde F_1,\tilde G_1):$ $C^{4,\alpha}(\bar\Omega)$ $\rightarrow$ $C^{1,\alpha}(\bar\Omega) \times C^{2,\alpha}(\partial\Omega) \times C^{2,\alpha}(\partial\Omega)$, $L_t^{-1}\circ G$ is a legitimate homotopy for the Leray-Schauder degree. It follows that
\begin{equation}
d = \deg_{L.S.}(L^{-1}\circ (\tilde F_1,\tilde G_1),\mathcal{\tilde O},0) = \deg_{L.S.}(L_1^{-1}\circ (\tilde F_1,\tilde G_1),\mathcal{\tilde O},0).
\label{LL::Eqn02}
\end{equation}
Next, set
\[
\tilde L_t\,w = \Big(\Delta F_3[w] - (1-t)\,N\,F_3[w], \Big(\frac{\partial \tilde F[w]}{\partial \gamma} + t\,\tilde F[w]\Big)\Big|_{\partial\Omega}, \Big(\Delta_T G_3[w] - ((1-t)N + t) G_3[w]\Big)\Big|_{\partial\Omega} \Big).
\]
Arguing as before, we have $\tilde L_t$ is invertible and 
\begin{equation}
\deg_{L.S.}(L_1^{-1}\circ (\tilde F_1,\tilde G_1),\mathcal{\tilde O},0) = \deg_{L.S.}(\tilde L_0^{-1}\circ (\tilde F_1,\tilde G_1),\mathcal{\tilde O},0) = \deg_{L.S.}(\tilde L_1^{-1}\circ (\tilde F_1,\tilde G_1),\mathcal{\tilde O},0).
\label{LL::Eqn03}
\end{equation}

Note that $\tilde L_1$ $=$ $(S\circ F_3, T \circ G_3)$ and so $\tilde L_1^{-1} \circ (\tilde F_1,\tilde G_1)$ $=$ $(F_3, G_3)^{-1} \circ (F_1, G_1)$. Hence, by \eqref{LL::Eqn02} and \eqref{LL::Eqn03}
\[
d = \deg_{L.S.}((F_3, G_3)^{-1} \circ (F_1, G_1), \mathcal{O},0)
	= \deg_{L.S.}((F_1,G_1)^{-1} \circ (F_3, G_3), \mathcal{O},0).
\]
Set
\[
A_t = (F_1,G_1)^{-1} \circ \big[(1-t)(F_3, G_3) - t(Id,0)\big],
\]
where $(Id,0)$ is considered as an operator from $C^{4,\alpha}(\bar\Omega)$ into $C^{2,\alpha}(\bar\Omega) \times C^{1,\alpha}(\partial\Omega)$. By the maximum principle and obliqueness, $A_t$ is a continuous family of invertible linear operators acting on $C^{4,\alpha}(\bar\Omega)$. Moreover, for $t$ $\in$ $[0,1)$, $(1-t)^{-1}\,A_t$ is of the form $\textrm{Id} + \textrm{Compact}$. Hence, by the homotopy invariance property of the Leray-Schauder degree,
\[
d = \deg((1-t)^{-1}A_t,\mathcal{\tilde O},0) \text{ for any } t \in [0,1),
\]
which implies
\[
d = (-1)^{\dim E^-(A_t)} \text{ for any } t \in [0,1),
\]
where
\[
E^-(A_t) = \bigoplus_{\lambda_i(t) < 0} \Big\{u \in C^{4,\alpha}(\bar\Omega): A_t\,u = \lambda_i(t)\,u\Big\}.
\]

To proceed, we claim that there exists some $C$ $>$ $0$ and $\delta$ $\in$ $(0,1)$ such that, for any $t$ $\in$ $(1-\delta,1]$
\begin{equation}
-C < \lambda < -\frac{1}{C} \text{ for any negative eigenvalue $\lambda$ of $A_t$}.
\label{EigenBnd}
\end{equation}
Indeed, let $\lambda$ be an eigenvalue of some $A_t$ and $u$ be a corresponding eigenfunction. Since $A_t$ is invertible, 
\[
\left\{\begin{array}{l}
a_{ij}\,u_{ij} + b_i\,u_i + c\,u = \frac{1}{\lambda}\big[(1-t)a_{ij}\,u_{ij} - tu\big] \text{ in } \Omega,\\
\beta_i\,u_i + \gamma\,u = \frac{1}{\lambda}\,(1-t)(\beta_i\,u_i + u) \text{ on } \partial\Omega,
\end{array}\right.
\]
which is equivalent to
\[
\left\{\begin{array}{l}
a_{ij}\,u_{ij} + \frac{\lambda}{\lambda - (1-t)}\,b_i\,u_i + \frac{\lambda}{\lambda - (1-t)}\,c\,u + \frac{t}{\lambda - (1-t)}u = 0 \text{ in } \Omega,\\
\beta_i\,u_i + \frac{\lambda}{\lambda - (1-t)}\gamma\,u - \frac{1-t}{\lambda - (1-t)}\,u = \text{ on } \partial\Omega,
\end{array}\right.
\]
It is readily seen that the first inequality in \eqref{EigenBnd} follows from the invertability of $(F_1,G_1)$ while the second follows from Lemma \ref{SemiFiniteness} for $\delta$ sufficiently small.

By \eqref{EigenBnd} and the compactness of $A_1$, we can pick a (simply connected) neighborhood $\mathcal{N}$ of $[-C,-\frac{1}{C}]$ in the complex plane such that in the set of eigenvalues of $A_1$ lying in $\mathcal{N}$ consists of all negative real eigenvalues of $A_1$. Furthermore, we can assume that $\mathcal{N}$ is symmetric about the real axis. Set
\begin{align*}
E(A_t,\mathcal{N})
	&= \bigoplus_{\lambda_i(t) \in \mathcal{N}}\Big\{u \in C^{4,\alpha}(\bar\Omega): A_t\,u = \lambda_i(t)u\Big\},\\
E^*(A_t,\mathcal{N})
	&= \bigoplus_{\lambda_i(t) \in \mathcal{N}\setminus \RR}\Big\{u \in C^{4,\alpha}(\bar\Omega): A_t\,u = \lambda_i(t)u\Big\}
\end{align*}
By the continuity of a finite system of eigenvalues (see e.g. \cite[pp. 213-214]{Kato}), for $\delta$ $>$ $0$ sufficiently small,
\[
\dim E(A_t,\mathcal{N}) \text{ is independent of } t \in (1-\delta,1].
\]
Also, since $A_t$ has real coefficients,
\[
\dim E^*(A_t,\mathcal{N}) \text{ is even}.
\]
Therefore, by \eqref{EigenBnd},
\[
d = \lim_{t \rightarrow 1}(-1)^{\dim E^-(A_t)} = \lim_{t \rightarrow 1}(-1)^{\dim E(A_t,\mathcal{N})} = (-1)^{\dim E(A_1,\mathcal{N})} = (-1)^{\dim E^-(A_1)}.
\]
As $A_1$ $=$ $-(F_1,G_1)^{-1}\circ(Id,0)$, \eqref{LL::Eqn01} follows. The proof of Corollary \ref{Cor:LSExt} is complete.
\end{proof}

Finally, we prove the uniqueness of the degree under properties $(p1)$--$(p3)$. We will only sketch the argument since it is standard. Let $d((F,G),\mathcal{O},0)$ be a degree which satisfies $(p1)$--$(p3)$. We will show that $d((F,G),\mathcal{O},0) = \deg((F,G),\mathcal{O},0)$. First, by Smale's infinite dimensional version of Sard's theorem, there exists $f_0 \in C^{2,\alpha}(\bar\Omega)$ and $g_0 \in C^{3,\alpha}(\partial\Omega)$ such that all zeroes of $(F - f_0,G - g_0)$ are non-degenerate and $\|F[u] - tf_0\|_{C^{2,\alpha}(\bar\Omega)} + \|G[u] - tg_0\|_{C^{3,\alpha}(\partial\Omega)} \geq c_0 > 0$ for all $u \in \partial\mathcal{O}$, all $t \in (0,1)$ and some $c_0 > 0$. By the homotopy invariance property $(p2)$, $d((F,G),\mathcal{O},0) = d((F - f_0,G - g_0),\mathcal{O},0)$ and $\deg((F,G),\mathcal{O},0) = \deg((F - f_0,G - g_0),\mathcal{O},0)$. Thus, we may assume from the beginning that all the zeroes of $(F,G)$ in $\mathcal{O}$ are non-degenerate. The uniqueness then follows from the addition property $(p1)$, Corollary \ref{Cor:a-d}(d), and the degree counting formula $(p3)$ for linear operators. 

We have finished the proof of Theorem \ref{Thm:DegExistence}. \hfill\qed

\section{Proof of Theorem \ref{Thm:LNBijection}}\label{Sec:Thm2}

We start with the injectivity of $L^N$:

\begin{proposition}\label{Prop:LNInjection}
Under the hypothesis of Theorem \ref{Thm:LNBijection}, there exists some constant $N_0$, depending only on $\|a_{st}\|_{C^{1,\alpha}}, \|b_i\|_{C^{1,\alpha}}, n,\lambda,\chi$ such that for all $N>N_0$, $L^N$ is injective.
\end{proposition}

\begin{proof} \noindent\underline{Step 1.} If $a_{st}w_{st} = 0$ in $\Omega$ and $L^N_{(3)} w = 0$ on $\partial\Omega$ simultaneously for some $N>0$ sufficiently large, then $V:= b_i\,D_i w + w = 0$ on $\partial \Omega$, hence $w \equiv 0$. 

Integrating by parts, we have
	\begin{equation}
	\begin{split}
		0 =& -\int_{\partial\Omega} \left(L^N_{(3)}w\right) V \\
		  =& -\int_{\partial\Omega} \left(b_i \Delta_Tw_i\right)\Big(\sum_j b_jD_jw + w\Big) + N\int_{\partial\Omega} |V|^2 \\
		\geq& \int_{\partial\Omega} |\sum_i b_i \nabla_T w_i|^2+N\int_{\partial\Omega}|V|^2 \\
		   &- \varepsilon\|\nabla_T Dw\|_{L^2(\partial\Omega)}^2-C_\varepsilon\left(\|Dw\|_{L^2(\partial\Omega)}^2+\|w\|_{L^2(\partial\Omega)}^2\right),
	\end{split}
	\label{M46}
	\end{equation}
for any positive constant $\varepsilon$, where the constant $C_\varepsilon>0$ depends on $\varepsilon$ and $\|b_i\|_{C^{1,\alpha}}$. Here $\nabla_T$ denotes the covariant gradient operator on $\partial\Omega$.

By Lemma \ref{Lem:Lem1}, we have
	\begin{equation*}
	\begin{split}
		&\int_{\partial\Omega} |Dw|^2 \leq C\int_{\partial\Omega} |V|^2,\\
		&\int_{\partial\Omega} |\nabla_T w_i|^2 \leq C\int_{\partial\Omega}\Big(\sum_j b_j \nabla_T w_j\Big)^2+  |V|^2.
	\end{split}
	\end{equation*}
Thus, by choosing $\varepsilon>0$ sufficiently small in \eqref{M46}, we deduce that
	\begin{equation}
		0\geq\frac12\int_{\partial\Omega}\sum_\tau\left(b_iD_{\tau}w_i\right)^2+\frac{N}{2}\int_{\partial\Omega}\left(b_iD_iw+w\right)^2.
	\end{equation}
This implies that $b_iD_iw+w=0$. Since $a_{st}D_{st} w = 0$ in $\Omega$, we obtain $w \equiv 0$ from the maximum principle. Step 1 is proved.

\bigskip
\noindent\underline{Step 2.} For any $w$ satisfying $L^N_{(3)}w = 0$, there holds
  	\begin{eqnarray}
		&& \|w\|_{H^2(\Omega)} \leq C\left(\|a_{st}w_{st}\|_{L^2(\Omega)} + \|w\|_{L^2(\Omega)}\right), \label{L3H2}\\
		&& \|w\|_{H^3(\Omega)} \leq C\left(\|a_{st}w_{st}\|_{H^1(\Omega)} + \|w\|_{L^2(\Omega)}\right),\label{L3H3}
	\end{eqnarray}
where $C$ depends on $\|a_{st}\|_{C^{1,\alpha}}$, $\|b_{i}\|_{C^{0,\alpha}}$, $n$, $\lambda$ and $\chi$.

Indeed, $L^N_{(3)}w=0$ implies that
  \begin{equation}\label{neop}
   L_T V:=-\triangle_T V + N V = g \quad\mbox{ on }\partial\Omega,
  \end{equation}
where $V = b_i\,w_i + w$ as in Step 1 and $g=(\triangle_T b_i)w_i+2\nabla_Tb_i\cdot\nabla_Tw_i+\triangle_T w$. The inverse operator $L_T^{-1} : H^{s-2}(\partial\Omega)\to H^{s}(\partial\Omega)$ is a bounded operator and $\|L_T^{-1}\| \leq C$ for some constant $C$ independent of $N$ (see Lemma \ref{Lem:BNorm} in the appendix). 
When $s=1/2$, by the trace theorem and interpolation we have
  \begin{equation}\label{H1/2}
  \begin{split}
   \|V\|_{H^{1/2}(\partial\Omega)} &\leq C\|g\|_{H^{-3/2}(\partial\Omega)} \\
	  &\leq C[\|w\|_{H^{1/2}(\partial\Omega)} + \|Dw\|_{H^{-1/2}(\partial\Omega)}]\\
	  &\leq C[\|w\|_{H^{1/2}(\partial\Omega)} + \|Dw\|_{H^{1/4}(\partial\Omega)}]\\
	  &\leq C\|w\|_{H^{7/4}(\Omega)} \leq \varepsilon \|w\|_{H^2(\Omega)} + C_\varepsilon\|w\|_{L^2(\Omega)},
  \end{split}
  \end{equation}
for any small $\varepsilon>0$.
Similarly, when $s=3/2$ we have
  \begin{equation}\label{H3/2}
   \|V\|_{H^{3/2}(\partial\Omega)} \leq \varepsilon \|w\|_{H^3(\Omega)} + C_\varepsilon\|w\|_{L^2(\Omega)}.
  \end{equation}

From $H^2, H^3$ estimates for linear elliptic equation of second order with oblique derivative boundary condition (see e.g. \cite[Theorem 15.2]{ADN}), we have
	\begin{eqnarray*}
		&& \|w\|_{H^2(\Omega)} \leq C\left(\|a_{st}w_{st}\|_{L^2(\Omega)} +\|V\|_{H^{1/2}(\partial\Omega)} + \|w\|_{L^2(\Omega)}\right), \\
		&& \|w\|_{H^3(\Omega)} \leq C\left(\|a_{st}w_{st}\|_{H^1(\Omega)}+\|V\|_{H^{3/2}(\partial\Omega)}  + \|w\|_{L^2(\Omega)}\right).
	\end{eqnarray*}
Recalling \eqref{H1/2} and \eqref{H3/2} we arrive at \eqref{L3H2} and \eqref{L3H3}.

\bigskip
\noindent\underline{Step 3.} For any $w$ satisfying $L^N_{(3)}w = 0$, there holds
  	\begin{eqnarray}
		&& \|w\|_{H^2(\Omega)} \leq C\|a_{st}w_{st}\|_{L^2(\Omega)}, \label{L3H2Strong}\\
		&& \|w\|_{H^3(\Omega)} \leq C\|a_{st}w_{st}\|_{H^1(\Omega)}.\label{L3H3Strong}
	\end{eqnarray}

We will only derive \eqref{L3H2Strong}. The proof of \eqref{L3H3Strong} is similar.

By \eqref{L3H2}, it suffices to show that 
\[
\|w\|_{L^2(\Omega)} \leq C\|a_{st}w_{st}\|_{L^2(\Omega)} \text{ for all } w \in Ker(L^N_{(3)}).
\]
Otherwise, there is a sequence $w^{(n)} \in Ker(L^N_{(3)})$ such that
\[
\|w^{(n)}\|_{L^2(\Omega)} = 1 \text{ but } \|a_{st} w^{(n)}_{st}\|_{L^2(\Omega)} \leq \frac{1}{n}.
\]
By \eqref{L3H2}, $\|w^{(n)}\|_{H^2(\Omega)} \leq C$, thus $w^{(n)}$ converges weakly in $H^2$ and strongly in $L^2$ to some $w^*$. It is straightforward to show that $a_{st}\,w^*_{st} = 0$ in $\Omega$ and $L^N_{(3)}w^* = 0$ on $\partial\Omega$. By Step 1, $w^* \equiv 0$. This contradicts $\|w^{(n)}\|_{L^2(\Omega)} = 1$ and the strong convergence of $w^{(n)}$ to $w^*$ in $L^2$. This concludes Step 3.

\bigskip
\noindent\underline{Step 4.} $L^N w = 0$ implies $w \equiv 0$.

Using $L^N_{(1)}w=0$, we have
	\begin{equation}
	\begin{split}
		0 &= \int_\Omega \left(L^N_{(1)}w\right)\left(-a_{kl}w_{kl}\right)\ dx \\
		  &\geq c_1\int_\Omega|D(a_{st}w_{st})|^2\ dx + N\int_\Omega|a_{st}w_{st}|^2\ dx - c_2\|w\|_{H^3(\Omega)}\|w\|_{H^2(\Omega)} \\
		  &\geq \frac{c_1}{2}\int_\Omega|D(a_{st}w_{st})|^2\ dx + (N-c_3)\int_\Omega|a_{st}w_{st}|^2\ dx,
	\end{split}
	\label{19}
	\end{equation}
where the constants $c_1,c_2,c_3$ depend only on $\|a_{st}\|_{C^{1,\alpha}}, \|b_i\|_{C^{1,\alpha}}$, $n$, $\lambda$ and $\chi$. In the above, the first inequality follows from integration by parts, $L^N_{(2)}w=0$, and H\"older's inequality, while the second inequality follows from the Cauchy-Schwarz inequality and \eqref{L3H2Strong}-\eqref{L3H3Strong}. It follows from \eqref{19} that for $N$ sufficiently large, we must have $a_{st}w_{st}\equiv0$. By Step 1, this implies $w \equiv 0$, which completes the proof.
\end{proof}

We now turn to proving the surjectivity of $L^N$. Consider $L_0^N : C^{4,\alpha}(\overline\Omega) \to C^\alpha(\overline\Omega) \times C^{1,\alpha}(\partial\Omega) \times C^{1,\alpha}(\partial\Omega)$ defined by 
	\begin{equation*}
		L_0^N w = \left(\triangle^2w - N \Delta w,\ \left.\left(\gamma_i D_i \Delta w)\right)\right|_{\partial\Omega},\ \left.\left(\gamma_i \triangle_T\left(D_i w\right)- N\gamma_i\,D_i w - Nw\right)\right|_{\partial\Omega}\right).
	\end{equation*}
Note that $(1-t)L_0^N + L^N$ coincides with $L^N$ when $a_{st}$ and $b_i$ are replaced with $(1-t)\delta_{st} + ta_{st}$ and $(1-t) \gamma_i + tb_i$, respectively. Therefore $(1-t)L_0^N + L^N$ is injective for all $t \in [0,1]$. By the continuity method \cite[Theorem 12.5]{ADN}, $L^N$ is isomorphic if and only if $L_0^N$ is isomorphic. In other words, it suffices to establish the surjectivity of $L_0^N$.

Define $\tilde L_0^N : C^{4,\alpha}(\overline\Omega) \to C^\alpha(\overline\Omega) \times C^{1,\alpha}(\partial\Omega) \times C^{1,\alpha}(\partial\Omega)$ defined by 
	\begin{equation*}
		\tilde L_0^N w = \left(\triangle^2w - N \Delta w,\ \left.\left(\gamma_i D_i \Delta w)\right)\right|_{\partial\Omega},\ \left.\left(\triangle_T\left(\gamma_i D_i w + w\right)- N(\gamma_i\,D_i w + w)\right)\right|_{\partial\Omega}\right),
	\end{equation*}
and let $\tilde L_t^N = (1 - t)\tilde L_0^N + t L_0^N$. It is easy to see that $\tilde L_0^N$ is bijective. By the continuity method \cite[Theorem 12.5]{ADN}, in order to show that $L_0^N$ is bijective, it suffices to show that $\tilde L_t^N$ is injective for all $0 \leq t \leq 1$. This is a consequence of the following lemma:

\begin{lemma}\label{Lem:Step1Ex}
If $\Delta w = 0$ in $\Omega$ and if
\[
M_tw := (1-t) \triangle_T\left(\gamma_i D_i w + w\right) + t\,\gamma_i \triangle_T\left(D_i w\right) - N(\gamma_i\,D_i w + w) = 0 \text{ on } \partial \Omega
\]
 for some $0 \leq t \leq 1$, then $V:= \gamma_i\,D_i w + w = 0$ on $\partial \Omega$, hence $w \equiv 0$. 

\end{lemma}

To see this, we adapt the argument in Step 1 of the proof of Proposition \ref{Prop:LNInjection}. We compute
	\begin{equation}
	\begin{split}
		0 =& -\int_{\partial\Omega} M_t w\, V \\
		  =& (1 - t) \int_{\partial_\Omega} |\nabla_T V|^2 - t\int_{\partial\Omega} \left(\gamma _i \Delta_T D_i w\right)V + N\int_{\partial\Omega} |V|^2 \\
		\geq& \int_{\partial\Omega} |\nabla_T V|^2+N\int_{\partial\Omega}|V|^2 \\
		   &- \varepsilon\|\nabla_T Dw\|_{L^2(\partial\Omega)}^2-C_\varepsilon\left(\|Dw\|_{L^2(\partial\Omega)}^2+\|w\|_{L^2(\partial\Omega)}^2\right),
	\end{split}
	\label{M46}
	\end{equation}
for any positive constant $\varepsilon$, where the constant $C_\varepsilon>0$ depends on $\varepsilon$ and $\|\gamma_i\|_{C^{1,\alpha}}$.
By Lemma \ref{Lem:Lem1}, we have
	\begin{equation*}
	\begin{split}
		&\int_{\partial\Omega} |Dw|^2 \leq C\int_{\partial\Omega} |V|^2,\\
		&\int_{\partial\Omega} |\nabla_T D w|^2 \leq C\int_{\partial\Omega} |\nabla_T V|^2+  |V|^2.
	\end{split}
	\end{equation*}
Thus, by choosing $\varepsilon>0$ sufficiently small in \eqref{M46}, we deduce that
	\begin{equation}
		0\geq\frac12\int_{\partial\Omega} |\nabla_T V|^2 +\frac{N}{2}\int_{\partial\Omega} |V|^2,
	\end{equation}
which implies that $V=0$ as desired. This concludes the proof of the lemma and hence of Theorem \ref{Thm:LNBijection}.\hfill\qed

%+++++++++++++++%

\section{Some applications}\label{Sec3}

\subsection{Boundary Yamabe problem}

Let $(M,g)$ be a smooth Riemannian manifold of dimension $n \geq 3$ and with non-empty boundary $\partial M$. Define the Schouten tensor
\[
A_g = \frac{1}{n - 2}\Big(\Ric_g - \frac{1}{2(n-1)}R_g\,g\Big),
\]
where $\Ric_g$ and $R_g$ are respectively the Ricci curvature and the scalar curvature of $g$. Let $h_g$ denote the mean curvature of $\partial M$. In conformal geometry, one is interested in finding a positive function $u$ such that the metric $g_u := u^{\frac{4}{n-2}}g$ such that 
\begin{equation}
\left\{\begin{array}{l}
f(\lambda(A_{g_u})) = 1, \lambda(A_{g_{u}}) \in \Gamma,\\
h_{g_u} = c,
\end{array}\right.
	\label{FNLYam}
\end{equation}
where $c$ is a given constant in $\RR$, $\Gamma \subset \RR^n$ satisfies
\begin{align*}
&\text{$\Gamma$ is an open, convex, symmetric cone with vertex at the origin,}\\
&\Gamma_1 := \{\lambda \in \RR^n: \sum\lambda_i > 0\} \supset \Gamma \supset \Gamma_n := \{\lambda \in \RR^n: \lambda_i > 0\},
\end{align*}
and $f \in C^\infty(\Gamma) \cap C^0(\bar\Gamma)$ satisfies
\begin{align*}
&f > 0 \text{ in } \Gamma \text{ and } f = 0 \text{ on } \partial \Gamma,\\
&f_{\lambda_i} > 0 \text{ in } \Gamma.
\end{align*}
We refer readers to \cite{LL} for the literature on this problem.

In \cite{LL}, an existence theorem for \eqref{FNLYam} was established using a degree theory which is less general than the one considered in the present paper. We outline some of the arguments here.

First, we restrict ourselves to the case where
\[
\text{$g$ is locally conformally flat,}
\]
and
\[
\text{$\partial M$ is umbilic.}
\]
Under the assumption that 
\[
\lambda(A_g) \in \Gamma,
\]
it was shown in \cite{LL} that all (positive) solutions of \eqref{FNLYam}, if exist, satisfy
\[
|\ln u| \leq C(M,g,f,\Gamma,c).
\]
Under an additinal assumption that
\[
\text{$f$ is concave and $c \geq 0$},
\]
one then appeals to known local $C^1$ and $C^2$ estimates, and Evans-Krylov's theory to conclude that all solutions of \eqref{FNLYam} satisfy
\[
\|\ln u\|_{C^{4,\alpha}(M)} \leq C(M,g,f,\Gamma,c,\alpha).
\]
It should be noted that, when $c < 0$, $C^2$ estimates fail at both local and global levels. See \cite{LL,LLCounter}.

With the above a priori estimate, our degree is defined and independent along a homotopy connecting $(f,\Gamma)$ to $(\sigma_1,\Gamma_1)$. Here $\sigma_1$ is the first elementary symmetric function. By property (e), the degree of $(\sigma_1,\Gamma_1)$ is the same as the Leray-Schauder degree for $(\sigma_1,\Gamma_1)$, which was computed to be non-zero in \cite{HanLi}. The desired existence result is established.

\subsection{Near-field reflector problem}

Consider the Monge-Amp\`ere type equation arising in a near-field reflector problem \cite{KW,LT}
	\begin{equation}\label{r01}
		\rho^*(T_u)\ \det\,DT_u=\rho(\cdot)\quad\mbox{ in }\Omega
	\end{equation}
with the boundary condition
	\begin{equation}\label{r02}
		T_u(\Omega)=\Omega^*,
	\end{equation}
where $\rho, \rho^*$ are, respectively, the intensities of incident and reflected rays satisfying
	\begin{equation}\label{r902}
		\int_\Omega\rho=\int_{\Omega^*}\rho^*,
	\end{equation}
$\Omega, \Omega^*$ are two bounded domains in $\mathbb{R}^n$ with $\overline\Omega\Subset B_1(0)$, and $T_u$ is the reflection mapping
	\begin{equation}\label{r03}
		T_u(x)=\frac{2Du}{|Du|^2-(u-Du\cdot x)^2},\quad x\in\Omega.
	\end{equation}

In \cite{LT} we proved the existence of solution for \eqref{r01} and \eqref{r02} under some convexity assumptions on domains $\Omega$ and $\Omega^*$ by using the degree theory established in Section 2. We outline the main steps here and refer the readers to \cite{LT} for more details.

First we show the a priori estimate that under suitable convexity and smoothness assumptions, the boundary condition \eqref{r02} is strictly oblique, and moreover, for any classical solution $u$ of \eqref{r01}--\eqref{r02}, we have the a priori estimate $\|u\|_{C^{4,\alpha}(\overline\Omega)}\leq C$ for $\alpha\in(0,1)$ and a positive constant $C$ depending on the given data.

To construct the homotopy family, we use the domain foliation, namely there exist an increasing family of smooth domains $\{\Omega_t\}, \{\Omega^*_t\}$, such that $\Omega_0=B_r(0)$, $\Omega^*_0=T_{u_0}(\Omega_0)$; $\Omega_1=\Omega$, $\Omega^*_1=\Omega^*$; and $\Omega_t, \Omega^*_t$ are uniformly convex and $Y^*$-convex, respectively, for all $0\leq t\leq 1$.

Consider the family of problems
	\begin{eqnarray}\label{r10}
		&& F_t[u_t] = \det\,[DT_{u_t}]-e^{\varepsilon(u_t-u_0)}\left\{\frac{t\rho}{\rho^*(T_{u_t})}+(1-t)\,\det\,[DT_{u_0}]\right\} =0\quad\mbox{in }\Omega_t, \\
		&& G_t[u_t] = \varphi^*_t\circ T(\cdot,u_t,Du_t)=0\quad\mbox{on }\partial\Omega_t, \nonumber
	\end{eqnarray}
where $\varepsilon>0$ is small, and $\varphi^*_t$ is the defining function of $\Omega^*_t$. 
Let $\varPhi:=\{\Phi_t : \mathbb{R}^n\to\mathbb{R}^n\}$ be a family of diffeomorphisms such that $\Phi_t(\Omega_t)=B_1(0)$ for each $t\in[0,1]$ and $\Phi_t\in C^5$ uniformly with respect to $t$. Define
	\begin{eqnarray}\label{r910}
		&& \tilde F_t[u] = F_t[u\circ\Phi_t]\quad\mbox{in }B_1(0), \\
		&& \tilde G_t[u] = G_t[u\circ\Phi_t]\quad\mbox{on }\partial B_1(0), \nonumber
	\end{eqnarray}
for any $u\in C^{4,\alpha}(\overline B_1)$.
 
Let $\mathcal{O}$ be a bounded open set in $\{u\in C^{4,\alpha}(\overline B_1)\,:\,\|u\|_{C^{4,\alpha}(\overline B_1)}\leq C(C\varepsilon+1)\}$ such that $\tilde F_t$ is elliptic and $\tilde G_t$ is oblique on $\mathcal{O}$, and $\partial\mathcal{O}\cap(\tilde F_t,\tilde G_t)^{-1}(0)=\emptyset$ for all $t\in[0,1]$. From the initial construction, $u_0$ is the unique solution of $(F_0,G_0)[u]=0$.
By the properties of the degree, $\deg((\tilde F_t,\tilde G_t),\mathcal{O},0)$ is defined for $0\leq t\leq 1$ and
	\begin{equation}\label{r11}
	\begin{split}
		\deg((\tilde F_1,\tilde G_1),\mathcal{O},0) &= \deg((\tilde F_0,\tilde G_0),\mathcal{O},0) \\
			&= \deg((F_0,G_0),\mathcal{O}_0,0)\neq0,
	\end{split}
	\end{equation}
where $\mathcal{O}_0=\{u\in C^{4,\alpha}(\overline\Omega_0)\,:\,\|u\|_{C^{4,\alpha}(\overline\Omega_0)}\leq C_\varepsilon+1\}$. This implies that there exists a solution $\tilde u_\varepsilon\in C^{4,\alpha}(\overline B_1)$ of the boundary value problem \eqref{r910} at $t=1$.
Hence there exists a solution $u_\varepsilon\in C^{4,\alpha}(\overline\Omega)$ of the boundary value problem
	\begin{eqnarray}\label{r12}
		&& \det\,[DT_u] = e^{\varepsilon(u-u_0)}\frac{\rho}{\rho^*(T_u)}, \\
		&& T_u(\Omega) = \Omega^* \nonumber
	\end{eqnarray}
for arbitrary small $\varepsilon>0$. To complete the existence proof we now need to let $\varepsilon\to0$. 
Write equation \eqref{r12} in the form of
	\begin{equation}\label{r13}
		\rho^*(T_u)\ \det\,DT_u = e^{\varepsilon(u-u_0)}\rho(\cdot)\quad\mbox{ in }\Omega.
	\end{equation}
Let $\{u_\varepsilon\}$ be the family of solutions of the problems \eqref{r13}. From \eqref{r902}
	\begin{equation*}
		\int_\Omega\rho=\int_{\Omega^*}\rho^*=\int_\Omega e^{\varepsilon (u_\varepsilon-u_0)}\rho,
	\end{equation*} 
we see that $u_\varepsilon-u_0$ must be zero somewhere in $\Omega$. Hence, from the a priori estimates $\|u_\varepsilon\|_{C^{4,\alpha}(\overline\Omega)}$ is bounded independently of $\varepsilon$. Thus a subsequence of $\{u_\varepsilon\}$ converges in $C^{4,\beta}(\overline\Omega)$ for $0<\beta<\alpha$ to a solution $u$ solving \eqref{r01}--\eqref{r02}, as required.

%%%%%%%%%%%%%%%%%5

\appendix
\section{Some well-known properties of the Laplace operator on a compact manifold}\label{AppA}

Let $(M,g)$ be a compact Riemannian manifold and $\Delta$ denote the Laplace operator of $g$. Let $H^s = H^s(M)$ denote the Sobolev space of some real exponent $s$. It is well-known that $-\Delta$ maps $H^s$ into $H^{s-2}$. In this appendix, we collect some well-known properties on the spectrum of $-\Delta: H^s \rightarrow H^{s-2}$ for $s \in [0,2]$ which are needed in the body of the paper.

For $s = 0$, $H^0 = L^2$. For $s = 2$, $H^2$ is defined as the set of space of $L^2$ functions whose first and second derivatives also belong to $L^2$. $H^{-2}$ is defined as the dual of $H^2$. For $s \in (0,2)$, $H^s$ can be defined as an interpolation space of $H^0 = L^2$ and $H^2$ as follows (see e.g. \cite{BerghLofstrom}):
\[
H^s = \Big\{u = \sum_{i = -\infty}^\infty u_i \in L^2 \Big| u_i \in H^2, \sum_{i=-\infty}^\infty (2^{-is}\|u_i\|_{L^2}^2 + 2^{i(1 - s)}\|u_i\|_{H^2}^2) < \infty\Big\}, \qquad 0 < s < 2.
\]
For $u \in H^s$, its $H^s$-norm is defined by
\[
\|u\|_{H^s}^2 = \inf\sum_{i=-\infty}^\infty (2^{-is}\|u_i\|_{L^2}^2 + 2^{i(2 - s)}\|u_i\|_{H^2}^2),
\]
where the infimum is taken over all possible represenation $u = \sum u_i$. For $s \in (-2,0)$, $H^s$ is defined as the dual of $H^{|s|}$, which is the same as the following interpolation space of $H^{-2}$ and $L^2$:
\[
H^s = \Big\{u = \sum_{i = -\infty}^\infty u_i \in H^{-2} \Big| u_i \in L^2, \sum_{i=-\infty}^\infty (2^{-i(2-|s|)}\|u_i\|_{H^{-2}}^2 + 2^{i|s|}\|u_i\|_{L^2}^2) < \infty\Big\}, \qquad -2 < s < 0.
\]
Its norm is defined similarly.

For $u \in H^2$, $\Delta u \in L^2$ is defined by standard differentiation. For $u \in L^2$, $\Delta u \in H^{-2}$ is defined by
\[
\left<\Delta u, v\right>_{H^{-2},H^2} = \left<u,\Delta v\right>_{L^2,L^2} \text{ for all } v \in H^2.
\]
Clearly $\Delta: H^2 \rightarrow L^2$ and $\Delta: L^2 \rightarrow H^{-2}$ are bounded linear operator. Furthermore, it is easy to check that $\Delta$ maps $H^s$ into $H^{s-2}$ for $s \in (0,2)$ and is a bounded linear operator between these spaces. Furthermore, for $\varphi \in C^\infty$,
\[
\left<\Delta u, \varphi \right>_{H^{s-2},H^{2 - s}} = \left<u, \Delta\varphi\right>_{L^2,L^2} \text{ for all } \varphi \in C^\infty(M).
\]

Next, fix some $N \geq 1$ and $s \in [0,2]$. Consider the operator $-\Delta + N: H^s \rightarrow H^{s-2}$. We claim that this operator is an isomorphism. This is clear for $s = 2$. For $s < 2$, note that $-\Delta + N: H^s \rightarrow H^{s-2}$ is injective. Indeed, if $-\Delta u + Nu = 0$ for some $u \in H^s$, then, for any $\psi \in C^\infty(M)$, there exists $\varphi \in C^\infty$ such that $(-\Delta + N)\varphi = \psi$ by standard elliptic regularity, and so
\[
\left<u,\psi\right>_{L^2,L^2} = \left<u,(-\Delta + N)\varphi\right>_{L^2,L^2} = \left<(-\Delta + N)u, \varphi\right>_{H^{s-2},H^{2-s}} = 0
\]
which implies that $u \equiv 0$. But as the formal adjoint of $-\Delta + N:H^s \rightarrow H^{s-2}$ is $-\Delta + N: H^{2 -s} \rightarrow H^{-s}$, which is injective, the Fredholm alternative implies that $-\Delta + N: H^s \rightarrow H^{s-2}$ is surjective. The claim is proved.

\begin{lemma}\label{Lem:BNorm}
There exists some constant $C$ such that for any $N \geq 1$ and $s \in [0,2]$, there holds
\[
\|(-\Delta + N)^{-1}\|_{\mathcal{L}(H^{s-2},H^s)} \leq C.
\]
\end{lemma}

\begin{proof}
Consider first $s = 2$. If $(-\Delta + N) v = u \in L^2$, then by integrating by parts and Cauchy-Schwarz's inequality
\[
\|D v\|_{L^2}^2 + N\|v|_{L^2}^2 = \left<u,v\right>_{L^2,L^2} \leq \frac{1}{2N}\|u\|_{L^2}^2 + \frac{N}{2}\|v\|_{L^2}^2,
\]
which implies
\[
\|D v\|_{L^2}^2 + N\|v|_{L^2}^2 \leq \left<u,v\right>_{L^2,L^2} \leq \frac{C}{N}\|u\|_{L^2}^2.
\]
This together with standard elliptic theory leads to
\[
\|D^2 v\|_{L^2}^2 \leq C(\|\Delta v\|_{L^2}^2 + \|v\|_{L^2}^2) = C(\|N v - u\|_{L^2}^2 + \|v\|_{L^2}^2) \leq C\|u\|_{L^2}^2.
\]
We have thus shown that
\[
\|(-\Delta + N)^{-1}\|_{\mathcal{L}(L^2,H^2)} \leq C.
\]
By duality, this gives
\[
\|(-\Delta + N)^{-1}\|_{\mathcal{L}(H^{-2},L^2)} \leq C.
\]

Now for a general $s \in (0,2)$ and $u \in H^{s-2}$, we write $u = \sum u_i \in H^{s-2}$ with $u_i \in L^2$. Then $(-\Delta + N)^{-1} u = \sum v_i \in H^s$ with $v_i = (-\Delta + N)^{-1} u_i \in H^2$. We then have
\begin{align*}
\|(-\Delta + N)^{-1} u\|_{H^s} 
	&\leq \sum (2^{-is}\|v_i\|_{L^2} + 2^{i(2-s)}\|v_i\|_{H^2})\\
	&\leq C \sum (2^{-is}\|u_i\|_{H^{-2}} + 2^{i(2-s)}\|u_i\|_{L^2}).
\end{align*}
Since this is true for all possible representations $u = \sum u_i$, we thus arrive at
\[
\|(-\Delta + N)^{-1} u\|_{H^s}  \leq C\|u\|_{H^{s-2}},
\]
which is exactly the assertion.
\end{proof}

\subsection*{Acknowledgment.} YanYan Li's research is partially supported by NSF-DMS-1501004. Jiakun Liu's research is partially supported by the Australian Research Council DE140101366.

\end{document}